\documentclass[a4paper,12pt]{article}
\usepackage{amsmath,amssymb,amsfonts,amsthm,eucal}
\usepackage{graphicx}
\usepackage{tikz}

\newtheorem{proposition}{Proposition}
\newtheorem{lemma}{Lemma}
\newtheorem{theorem}{Theorem}
\newtheorem{definition}{Definition}

\newtheorem{axiom}{Axiom}

\begin{document}
\date{\today}
\title{On proportionality in multi-issue problems with crossed claims}

\author{Rick K. Acosta-Vega\thanks{Faculty of Engineering, University of Magdalena, Colombia. \{racosta@unimagdalena.edu.co\}}
\and Encarnaci\'{o}n Algaba\thanks{%
Department of Applied Mathematics II and IMUS, University of Sevilla, Spain. \{ealgaba@us.es\}}
\and Joaqu\'in S\'anchez-Soriano\thanks{%
\textbf{Corresponding author}. R.I. Center of Operations Research (CIO), Miguel Hern\'andez University of Elche, Spain. \{joaquin@umh.es\}}
}
\maketitle

\begin{abstract}
In this paper, we analyze the problem of how to adapt the concept of proportionality to situations where several perfectly divisible resources have to be allocated among certain set of agents that have exactly one claim which is used for all resources. In particular, we introduce the constrained proportional awards rule, which extend the classical proportional rule to these situations. Moreover, we provide an axiomatic characterization of this rule. 
\end{abstract}

\section{Introduction}\label{intro}

Allocation problems describe situations in which a resource (or resources) must be distributed among a set of agents. These problems are of great interest in many settings, for this reason the literature on the matter is extensive. A particular allocation problem is arisen in situations where there is a perfectly divisible resource over which there is a set of agents who have rights or demands, but the resource is not sufficient to satisfy them. This problem is known as bankruptcy problem and was first formally analyzed in O'Neill (1982) and Aumann and Maschler (1985). Since then it has been extensively studied in the literature and many allocation rules have been defined (see Thomson, 2003, 2015, 2019, for a detailed inventory of rules). In the literature many applications of bankruptcy problems can be found. Some examples are the following. Pulido et al. (2002, 2008) study allocation problems in university management; Niyato and Hossain (2006), Gozalvez et al. (2012), and Lucas-Estañ et al. (2012) analyze radio resource allocation problems in telecommunications; Casas-Mendez et al. (2011) study the musseum pass problem; Hu et al. (2012) analyze the airport problem; Giménez-Gómez et al. (2016), Gutiérrez et al. (2018), and Duro et al. (2020) analyze the CO2 allocation problem; Sanchez-Soriano et al. (2016) study the apportionment problem in proportional electoral systems; and Wickramage et al. (2020) analyze water allocation problems in rivers.

An extension of bankruptcy problems are multi-issue allocation problems (Calleja et al., 2005). These describe situations in which there is a (perfect divisible) resource which can be distributed between several issues, and a (finite) number of agents that have claims on each of those issues, such that the total claim is above the available resource. This problem is also solved by means of allocation rules and there are different ways to do it (see, for example, Calleja et al. (2005), Borm et al. (2005), and Izquierdo and Timoner (2016)). Ju et al. (2007), Moreno-Ternero (2009) and Berganti\~nos et al. (2010) study the proportional rule for multi-issue allocation problems.

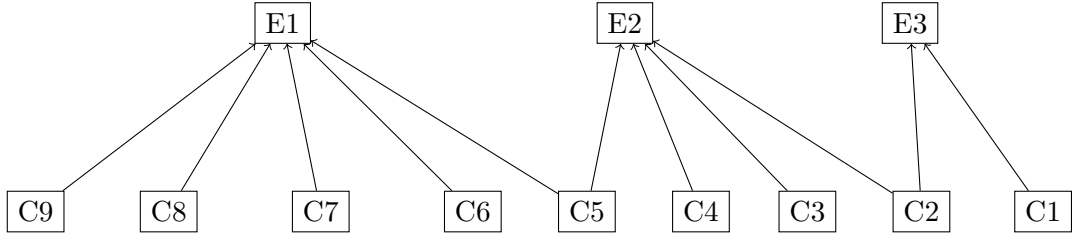
\begin{figure}
\begin{center}
\begin{tikzpicture}
\node[draw] (c8) at (10.5,0) {{\small C1}};
\node[draw] (c7) at (8.9,0) {{\small C2}};
\node[draw] (c6) at (7.4,0) {{\small C3}};
\node[draw] (c9) at (6,0) {{\small C4}};
\node[draw] (c5) at (4.5,0) {{\small C5}};
\node[draw] (c4) at (3,0) {{\small C6}};
\node[draw] (c3) at (1,0) {{\small C7}};
\node[draw] (c2) at (-1,0) {{\small C8}};
\node[draw] (c1) at (-2.75,0) {{\small C9}};

\node[draw] (e1) at (0.5,2.5) {{\small E1}};
\node[draw] (e2) at (5,2.5) {{\small E2}};
\node[draw] (e3) at (8.75,2.5) {{\small E3}};

\draw [->] (c1) -- (e1);
\draw [->] (c2) -- (e1);
\draw [->] (c3) -- (e1);
\draw [->] (c4) -- (e1);
\draw [->] (c5) -- (e1);

\draw [->] (c5) -- (e2);
\draw [->] (c6) -- (e2);
\draw [->] (c7) -- (e2);
\draw [->] (c9) -- (e2);

\draw [->] (c7) -- (e3);
\draw [->] (c8) -- (e3);

\end{tikzpicture}
\end{center}
\caption{Multi-issue problem with crossed claims.}\label{2levels}
\end{figure}

However, the situation described in Figure \ref{2levels} does not fit to a multi-issue allocation problem as referred in the previous paragraph, but to a multi-issue allocation problem with crossed claims introduced by Acosta-Vega et al. (2021). These describe situations in which there are several (perfect divisible) resources and a (finite) set of agents who have claims on them, but only one claim (not a claim for each resource) with which one or more resources are requested. The total claim for each resource exceeds its availability. This problem is solved by means of allocation rules (Acosta-Vega et al., 2021).

In this paper, in order to solve allocation problems as the described in Figure \ref{2levels}, we introduce the constrained proportional awards rule for multi-issue allocation problems with crossed claims that naturally extends the proportional rule for single issue allocation problems. This rule is characterized axiomatically by using five properties: \emph{Pareto efficiency}, \emph{equal treatment of equals}, \emph{guaranteed minimum award}, \emph{consistency},  and \emph{non-manipulability by splitting}. The first one says that there is no a feasible allocation in which at least one of the claimants receive more. \emph{Equal treatment of equals} states that equal agents must receive the same. \emph{Guaranteed minimum award} means that a claimant should not receive less than what she would receive in the worst case, if the issues were distributed separately. \emph{Consistency} requieres that if a subset of agents leave the problem respecting what had been allocated to those who remain, then what those agents receive in the new problem is the same as what they received in the original problem. Finally, \emph{non-manipulability by splitting} means that it is not profitable to split one agent in several agents. Moreover, although it is not necessary in characterization, the constrained proportional awards rule also satisfies \emph{non-manipulability by restricted merging} which guarantees that it is not profitable to merge several ``homologous'' agents into one.

The rest of the paper is organized as follows. Section \ref{model} presents multi-issue allocation problems with crossed claims (MAC), and the concept of rule for these problems. In Section \ref{prop}, the constrained proportional awards rule for multi-issue bankruptcy problems with crossed claims is defined. In Section \ref{properties}, we present several properties which are interesting in the context of MAC problems. In Section \ref{char}, we characterize the constrained proportional awards rule. Section \ref{conc} concludes.


\section{Multi-issue allocation problems with crossed claims}
\label{model}
A one-issue allocation problem is given by a triplet $(N,E,c)$, where $N$ is the set of claimants, $E \in \mathbb{R}_+$ is a perfectly divisible amount of resource (the \emph{issue} or \emph{estate}) to be divided, and $c$ is the vector of demands, such that $C = \sum_{j \in N} c_j> E$.  One of the most relevant ways to allocate the resource among the claimants in one-issue allocation problem is the proportional rule (PROP) (Aristotle, 4th Century BD), which is defined as follows:

\begin{equation}
\label{defPROP}
PROP_j(N,E,c)=\frac{c_j}{C}E,\quad j\in N.
\end{equation} 
 
We now consider a situation where there are a finite set of issues  $\mathcal{I}=\{1,2,...,m\}$ such that  there is a perfectly divisible amount $e_i$ of each issue $i$. Let $E =(e_1,e_2,...,e_m)$ be the vector of available amounts of issues. There are a finite set of claimants $N=\{1,2,...,n\}$ such that each claimant $j$ claims $c_j$. Let $c=(c_{1},c_{2},...,c_{n})$ be the vector of claims. However, each claimant claims to different sets of issues, in general. Thus, let $\alpha$ be a set-valued function that associates with every $j \in N$ a set $\alpha(j) \subset \mathcal{I}$. In fact, $\alpha(j)$ represents the issues to which claimant $j$ asks for. Furthermore, $\sum_{j:i \in \alpha(j)} c_j > e_i$, for all $i \in \mathcal{I}$, otherwise, those issues could be discarded from the problem because they do not impose any limitation, and so the allocation would be trivial. Therefore, a multi-issue allocation problem with crossed claims (MAC in short) is defined by a 5-tuple $(\mathcal{I},N,E,c,\alpha)$, and the family of all these problems is denoted by $\mathcal{MAC}$.

Given a problem $(\mathcal{I},N,E,c,\alpha) \in \mathcal{MAC}$, a \emph{feasible allocation} for it, it is a vector $x \in \mathbb{R}^N$ such that:
\begin{enumerate}
\item $0 \leq x_i \leq c_i$, for all $i \in N$.
\item $\sum_{i \in N:j \in \alpha(i)}x_i \leq e_j$, for all $j \in M$,
\end{enumerate}
and we denote by $A(\mathcal{I},N,E,c,\alpha)$ the set of all its feasible allocations.

These two requirements are standard in the literature of allocation problems.

A \emph{rule} for multi-issue bankruptcy problems with crossed claims is a mapping $R$ that associates with every $(\mathcal{I},N,E,c,\alpha) \in \mathcal{MBC}$ a unique feasible allocation $R(\mathcal{I},N,E,c,\alpha) \in A(\mathcal{I},N,E,c,\alpha)$.


\section{The constrained proportional awards rule for $\mathcal{MAC}$ problems}\label{prop}

The proportional rule (PROP) is perhaps the most important rule to solve allocation problems in general,\footnote{See, for instance, Algaba et al. (2019) who present two solutions belonging to the family of proportional solutions for the problem of sharing the profit of a combined ticket for a transport system.} and bankruptcy problems in particular. This rule simply divides the resource in proportion to the claims. The question in $\mathcal{MAC}$ problems is what ``in proportion to the claims'' means. In the context of one-issue allocation problems, ``in proportion to the claims'' means that all claimants receive the same amount for each unit of claim. How to extrapolate this to the MAC situations. To answer this question, we introduce the constrained proportional awards rule (CPA in short) as the result of an iterative process in which the available amount of at least one of the issues is fully distributed in each step and so on and so forth while possible. This rule is formally defined below.

\begin{definition}\label{CPA}
Let $(\mathcal{I},N,E,c,\alpha) \in \mathcal{MAC}$, the \emph{constrained proportional awards rule} for $(\mathcal{I},N,E,c,\alpha)$, $CPA(\mathcal{I},N,E,c,\alpha)$,  is defined by means of the following iterative procedure:
\begin{itemize}
\item[Step 0.]
\begin{enumerate}
\item $\mathcal{I}^1 = \{i \in \mathcal{I}: e_i^1>0\}$ is the set of active issues.
\item $\mathcal{N}^1 = \{j \in N: c_j^1 >0 \text{ and } e_i^1 >0, \forall i \in \alpha(j)\}$ is the set of active claimants.
\item For each $i \in \mathcal{I}$, $e_i^1 = e_i$, and for each $j \in N$, $c_j^1 = c_j$.
\end{enumerate} 
\item[Step s.] 
\begin{enumerate}
\item $\mathcal{N}^s = \{j \in N: c_j^s >0 \text{ and } e_i^s >0, \forall i \in \alpha(j)\}$. $\mathcal{I}^s = \{i \in \mathcal{I}: e_i^s >0\}$.
\item For each $i \in \mathcal{I}^s$, we calculate the greatest $\lambda_{i} ^s$, so that $\lambda_{i} ^s\sum_{j \in \mathcal{N}^s:i \in \alpha(j)}c_j^s \leq e_i^s$, and take $\lambda^s = \min\{\lambda_{i} ^s:i \in \mathcal{I}^s\}$. 
\item Now, we allocate to each claimant $j \in \mathcal{N}^s$, $a_j^s=\lambda^sc_j^s$, and $a_j^s =0$ to the non-active claimants.
\item We update the active issues, $\mathcal{I}^{s+1}$, and the active claimants, $\mathcal{N}^{s+1}$. If $\mathcal{I}^{s+1} = \varnothing$ or $\mathcal{N}^{s+1} = \varnothing$, then the process ends, and
$$
CPA_j(\mathcal{I},N,E,c,\alpha) = \sum_{h =1}^{s}a_j^h, \forall j\in N.
$$
Otherwise, the available amounts of issues and the claims are updated:
$$e_i^{s+1} = e_i^s - \lambda^s\sum_{j\in N:i \in \alpha(j)}c_j^s, \forall i \in \mathcal{I}, \text{ and }  \hspace{0.5em} c_j^{s+1} =c_j^s - \lambda^sc_j^s, \forall j \in N,
$$
and we go to Step $s+1$.
\end{enumerate}
\end{itemize}
\end{definition}
The iterative procedure of CPA is well-defined and always leads to a single point. Moreover, since in each step at least the available amount of one issue is distributed in its entirety, except maybe in the last step, it ends in a finite number of steps, at most $|\mathcal{I}|$. Finally, when we have a one-issue allocation problem, then we obtain PROP. Therefore, this definition extends PROP to the context of MAC.

From the application of the iterative process to calculate CPA, we can consider the chains of active issues and active claimants in the application of the procedure to calculate $CPA(\mathcal{I},N,E,c,\alpha)\in \mathcal{MAC}$:
$$
\mathcal{I}^1 \supset \mathcal{I}^2 \supset \ldots \supset \mathcal{I}^r, \text{ and } \mathcal{N}^1 \supset \mathcal{N}^2 \supset \ldots \supset \mathcal{N}^r
$$

From these chains, we can establish an order relationship between issues as follows. We say that issue $i_1$ strictly precedes issue $i_2$ in a chain of actives issues, $i_1 \prec i_2$, if there is $\mathcal{I}^s$ such that $i_1 \notin\mathcal{I}^s$ and $i_2 \in\mathcal{I}^s$, i.e., $i_1$ becomes non-active before than $i_2$. We write $i_1 \preceq i_2$ when $i_1$ becomes non-active before than $i_2$ or both issues become non-active at the same time. Finally, we write  $i_1 \simeq i_2$ when both issues become non-active at the same time. Analogously, we can establish an order relationship between claimants.

Furthermore, we can associate with each pair of sets $\mathcal{I}^s$ and $\mathcal{N}^s$ a number $\rho^s$, $\rho^s \in [0,1]$, which represents the proportion of claims obtained by claimants in $\mathcal{N}^s$ but not in $\mathcal{N}^{s+1}$. Moreover, by construction $\rho^s < \rho^{s+1}$. Thus, we have that
$$
0<\rho^1 < \rho^2 <\ldots< \rho^r\leq 1.
$$

These $\rho's$ represent the accumulative proportion of the claims allocated to the claimants, i.e., what part of their claims they have received up to a given step of the iterative procedure. In this way, this procedure is reminiscent of the constrained equal awards rule (CEA) in bankuptcy problems, but instead of using the principle of egalitarianism, the principle of proportionality is used, hence the name of \emph{constrained proportional awards} rule. Therefore, not all claimants receive the same proportion of their claims, but the rule tries to keep the proportionality as much as possible restricted to (1) the relation between the available amounts of issues and the total claims to them, and (2) the goal of allocating as much as possible of all available amounts of issues.

Below some results about the iterative process that defines CPA are given.

\begin{proposition}The following statements hold
\begin{enumerate}
\item Given $(\mathcal{I},N,E,c,\alpha) \in \mathcal{MAC}$, if there are problems $(\mathcal{I}_1,N_1,E_1,c^1,\alpha_1), (\mathcal{I}_2,N_2,E_2,c^2,\alpha_2) \in \mathcal{MAC}$, such that $\mathcal{I}_1 \cup \mathcal{I}_2 = \mathcal{I}$, $N_1 \cup N_2 = N$, $E_1 \oplus E_2 = E$,  $c^1 \oplus c^2 = c$, and $\alpha_1(j)=\alpha(j), \forall j \in N_1$ and $\alpha_2(j)=\alpha(j), \forall j \in N_2$, so that $\left(\bigcup_{j \in N_1}\alpha(j)\right) \bigcap \left(\bigcup_{j \in N_2}\alpha(j)\right) = \varnothing$, then
$$
CPA(\mathcal{I},N,E,c,\alpha) = CPA(\mathcal{I}_1,N_1,E_1,c^1,\alpha_1) \oplus CPA(\mathcal{I}_2,N_2,E_2,c^2,\alpha_2).\footnote{Given $a \in \mathbb{R}^n$ and $b \in \mathbb{R}^m$, $a \oplus b \in \mathbb{R}^{n+m}$, i.e., $\oplus$ is the concatenation operator of two vectors.}
$$
\item If $\lambda^s <1$, then each $h \in \arg\min\{\lambda_{i} ^s:i \in \mathcal{I}^s\} \subset \mathcal{I}$ becomes non-active in the next step.
\item If $\lambda^s = 1$ for some $s$, then the iterative procedure ends in that step.
\end{enumerate}
\begin{proof}
These statements follow from the own structure of the iterative procedure to calculate CPA. We next provide a brief outline of each one.
\begin{enumerate}
\item This follows from the fact that since there are no crossed demands between the two subproblems, they do not affect each other and, therefore, the results are independent of each other.
\item $\lambda^s < 1$ implies that there are some issues in $\mathcal{I}^s$ for which their available amounts are not enough to satisfy the total claims to them in step $s$. Thus, $\lambda_{h} ^s\sum_{j \in \mathcal{N}^s:h \in \alpha(j)}c_j^s = e_h^s,\hspace{0.5em} \forall h \in \arg\min\{\lambda_{i} ^s:i \in \mathcal{I}^s\}$. Therefore, when $\lambda^s$ is applied, all those issues become non-active in  the next step.
\item $\lambda^s = 1$ implies that all active claimants in $\mathcal{N}^{s}$ receive their pendent claims, then they become non-active in the next step, i.e., $\mathcal{N}^{s+1}=\varnothing$. Therefore, the procedure ends.
\end{enumerate}
\end{proof}
\end{proposition}

The first statement says that if a problem can be separated into two disjoint problems, then it is the same to calculate CPA for the whole problem as for each of them and then paste the results. The second states when an issue becomes non-active. Finally, the third provides another stopping criterion for the iterative procedure to calculate CPA. Moreover, when the procedure ends with $\lambda = 1$, then it means that there may be resources left over from some issues. Otherwise, all resources have been fully distributed.


\section{Properties}\label{properties}

In this section, we present several properties which are interesting in the context of MAC problems. These properties are related to efficiency, fairness, consistency, and manipulability.

First, we introduce two concepts related to two claimants comparisons. In MAC situations claimants are characterized by two elements: their claims and the issues to which they claim. Therefore, both should be taken into account when establishing comparisons among them.

\begin{definition}
Let $(\mathcal{I},N,E,c,\alpha) \in \mathcal{MAC}$, and two claimants $j,k \in N$, we say they are {\em homologous}, if $\alpha(j) = \alpha(k)$; and we say that they are {\em equal}, if they are homologous and $c_j=c_k$.
\end{definition}

Next, we give a set of properties which are very natural and reasonable for an allocation rule in MAC situations.

The first property relates to efficiency. In allocation problems is desirable that the resources to be fully distributed, but in MAC situations this is not always possible (see Acosta-Vega, 2021). Therefore, a weaker version of that is considered in which only is required that there is no a feasible allocation in which at least one of the claimants receive more. This is established in the following axiom.

\begin{axiom}[PEFF]\label{PEFF}
Given a rule $R$, it satisfies {\em Pareto efficiency}, if for every problem $(\mathcal{I},N,E,c,\alpha) \in \mathcal{MAC}$, there is no a feasible allocation $a \in \mathbb{R}^N_{+}$ such that $a_j \geq R_j(\mathcal{I},N,E,c,\alpha), \forall j \in N$, with at least one strict inequality.
\end{axiom}

Note that $PEFF$ implies that at least the available amount of one issue is fully distributed, and no amount is left of an issue undistributed, if it is possible to do so. However, it does not require that all available amounts of the issues have to be fully distributed. On the other hand, a feasible allocation that satisfies the condition in Axiom \ref{PEFF} is called Pareto efficient. 

The second property states that equal claimants should receive the same in the final allocation. This is a basic requirement of fairness and non-arbitrariness. This is defined in the following axiom.

\begin{axiom}[ETE]\label{ETE}
Given a rule $R$, it satisfies {\em equal treatment of equals}, if for every problem $(\mathcal{I},N,E,c,\alpha)\in \mathcal{MAC}$ and every pair of equal claimants $j,k \in N$, $R_j(\mathcal{I},N,E,c,\alpha)=R_k(\mathcal{I},N,E,c,\alpha)$.
\end{axiom}

The third property assures the minimum that should be guaranteed to each claimant. In our case, these minimum amounts are determined from the analysis of the problems associated with each issue independently. In particular, the property states that a claimant should not receive less than what she would have received in the worst case, if the rule had been applied to each problem separately to each of the issues. This is established in the following property.

\begin{axiom}[GMA]\label{GMA}
Given a rule $R$, it satisfies {\em guaranteed minimum award}, if for every problem $(\mathcal{I},N,E,c,\alpha)\in \mathcal{MAC}$,
$$
R_j(\mathcal{I},N,E,c,\alpha)\geq \min\left\{R_j\left(\{i\},N_i, e_i, c|_{N_i}\right) : i \in \alpha(j)\right\}, \forall j \in N,
$$
where $N_i = \{k \in N : i \in \alpha(k)\}$, and $c|_{N_i}$ is the vector whose coordinates correspond to the claimants in $N_i$.
\end{axiom}

The fourth property is a requirement of robustness when some agents leave the problem with their allocations (see Thomson, 2011, 2018). In particular, when a subset of claimants leave the problem respecting the allocations to those who remain, then it seems reasonable that claimants who leave will receive the same in the new problem as they did in the original. This is formally given in the following axiom.

\begin{axiom}[CONS]\label{CONS}
Given a rule $R$, it satisfies {\em consistency}, if for every problem $(\mathcal{I},N,E,c,\alpha)\in \mathcal{MAC}$, and $N' \subset N$, it holds that 
$$
R_j(\mathcal{I},N,E,c,\alpha) = R_j(\mathcal{I}',N',E'^R,c|_{N'},\alpha), \text{ for all } j \in N',
$$
where $(\mathcal{I}',N',E',c|_{N'},\alpha)\in \mathcal{MAC}$, called   the \emph{reduced problem associated with $N'$}, $\mathcal{I}'=\{i \in \mathcal{I}: \text{ there exists } k \in N' \text{such that }i \in \alpha(k)\}$, $E'^R=(e'^R_1,\ldots,e'^R_m)$ so that $e'^R_i=e_i- \sum_{j \in N\backslash N':i \in \alpha(j)}R_j(\mathcal{I},N,E,c,\alpha)$, for all $i \in \mathcal{I}'$, and $c|_{N'}$ is the vector whose coordinates correspond to the claimants in $N'$.
\end{axiom}

The last two properties are related to claimants' ability to manipulate the final allocation by splitting their claims among several new claimants or merging their claims into a single claimant. It seems sensible that if the claimants do this, they will not benefit and receive the same as they did in the original problem. These two possibilities are established in the following axioms.

\begin{axiom}[NMS]\label{NMS}
Given a rule $R$, it satisfies {\em non-manipulability by splitting}, if for every pair of problems $(\mathcal{I},N,E,c,\alpha), (\mathcal{I},N',E,c',\alpha') \in \mathcal{MAC}$, such that:
\begin{enumerate}
\item $N \subset N'$, $S = \{i_1, \ldots, i_k\}$, such that $N' = (N\backslash S) \cup S_{i_1} \cup \ldots \cup S_{i_m}$, where $S_{i_k}$ is the set of agents into which agent $i_k$ has been divided.  
\item $c'_j = c_j, \forall j \in N\backslash S$ and $\sum_{k \in S_{i_h}}c'_k =c_{i_h}, h =1,\ldots,m$,
\item $\alpha'(j) = \alpha(j), \forall j \in N\backslash S$ and $\alpha'(j) = \alpha(i_h), \forall j \in  S_{i_h}, h =1,\ldots,m$,
\end{enumerate}
it holds
$$
\sum_{j \in S_{i_h}}R_j(\mathcal{I},N',E,c',\alpha') =R_{{i_h}}(\mathcal{I},N,E,c,\alpha), h =1,\ldots,m .
$$
\end{axiom}

\begin{axiom}[NMRM]\label{NMRM}
Given a rule $R$, it satisfies {\em non-manipulability by restricted merging}, if for every pair of problems $(\mathcal{I},N,E,c,\alpha), (\mathcal{I},N',E,c',\alpha') \in \mathcal{MAC}$, such that:
\begin{enumerate}
\item $N \subset N'$,
\item $c_j = c'_j, \forall j \in N\backslash \{j_0\}$ and $c_{j_0} = \sum_{k \in (N'\backslash N) \cup \{j_0\}}c'_k$,
\item $\alpha(j) = \alpha'(j), \forall j \in N\backslash \{j_0\}$ and $\alpha(j) = \alpha'(j_0), \forall j \in (N'\backslash N) \cup \{j_0\}$,
\end{enumerate}
it holds
$$
R_{j_0}(\mathcal{I},N,E,c,\alpha) = \sum_{j \in (N'\backslash N) \cup \{j_0\}}R_j(\mathcal{I},N',E,c',\alpha').
$$
\end{axiom}

Note that in $NMS$ we move from the allocation problem with set of claimants $N$ to the problem with set of claimants $N'$, i.e., one of the claimants is splitted into several new claimants, one of whom has the same name as in $N$. However, in $NMRM$ we move from the problem in $N'$ to the problem in $N$, i.e., several claimants merge into one claimant who has the same name as in $N'$, but all merged claimants are homologous. Thus, we are only considering the merging of homologous claimants. For this reason we call this axiom non-manipulability by ``restricted'' merging. Nevertheless, it seems reasonable from a perspective of symmetry of both properties, because when one claimant is splitted into several new claimants, these are homologous in the new problem.

CPA satisfies all properties above mentioned. We establish this in the following theorem.

\begin{theorem}\label{CPAprop}
CPA for multi-issue bankruptcy problems with crossed claims satisfies $PEFF$, $ETE$, $GMA$, $CONS$, $NMS$, and $NMRM$.
\begin{proof}We go axiom by axiom.
\begin{itemize}
\item CPA satisfies $PEFF$ and $GMA$ by definition.
\item If two claimants are equal, then CPA allocates both the same, since the procedure to calculate the rule treats, in each step, all active equal claimants egalitarianly, so if two claimants are equal, they stop receiving at the same step. Therefore, CPA satisfies $ETE$.
\item Given $(\mathcal{I},N,E,c,\alpha)\in \mathcal{MAC}$ and $(\mathcal{I}',N',E'^{CPA},c|_{N'},\alpha)\in \mathcal{MAC}$ the reduced problem associated with $N' \subset N$. Let us consider the following sets obtained from the application of CPA to $(\mathcal{I},N,E,c,\alpha)$:
$$
\mathcal{A}^s=\mathcal{N}^s \setminus \mathcal{N}^{s+1}, s=1,\ldots,r, \text{ and } \mathcal{B}^s=\mathcal{I}^s \setminus \mathcal{I}^{s+1}, s=1,\ldots,r.
$$

We now consider the following sets: $N' \cap \mathcal{A}^s, s=1,\ldots,r$. Taking into account the definitions of $E'^{CPA}$, $\mathcal{N}^s$, and $\mathcal{I}^s$, it is evident that claimants in $N' \cap \mathcal{A}^s$ cannot receive more than $\rho^s$ times their claims, because otherwise, the available amounts of issues $e'^{CPA}_i$, $i \in \mathcal{B}^s$, would be exceeded. Thus, from the definition of CPA, claimants in $N' \cap \mathcal{A}^s$ have to receive exactly $\rho^s$ times their claims. Therefore, the claimants in $(\mathcal{I}',N',E'^{CPA},c|_{N'},\alpha)$ receive the same as in $(\mathcal{I},N,E,c,\alpha)$. Hence, CPA satisfies consistency.

\item Note that when one claimant splits into several new claimants, CPA for the new problem will have the same number of iterations as in the original one, since the claim for each issue will be obviously the same in each step. Therefore, all split claimants will receive the same proportion of their claims which coincides with the proportion obtained by the split claimant in the original problem. Thus, the aggregate allocation of the split claimants in the new problem coincides with the allocation of the split claimant in the original problem.
\item When two homologous claimants merge into a new one claimant, we can make a completely analogous reasoning as in the case of a claimant splits into several new claimants. Therefore, CPA also satisfies NMRM.
\end{itemize}
\end{proof}
\end{theorem}


\section{Characterization}\label{char}
In this section, the aim is to achieve a better knowledge of the CPA rule for $\mathcal{MAC}$ by describing it in a unique way as a combination of some reasonable axioms. We characterize the $CPA$ rule by means of $PEFF$, $ETE$, $GMA$, $CONS$, and $NMS$. Therefore, the CPA rule can be considered as a desirable way to distribute a set of issues among their claimants. Before giving the characterization of CPA, we need the following lemmas.

\begin{lemma}\label{lemma}
Let $(\mathcal{I},N,E,c,\alpha) \in \mathcal{MAC}$, such that $|\mathcal{I}|=1$. If a rule $R$ satisfies $PEFF$, $ETE$, and $NMS$, then 
$$
R_i(\mathcal{I},N,E,c,\alpha) = \frac{c_1}{\sum_{j \in N} c_j}e, \text{ for all } i \in N.
$$
\begin{proof} First note that in this case the function $\alpha$ is irrelevant. Let $R_1,R_2,\ldots, R_n$ be the allocations for claimants in $N$, respectively. By $PEFF$, and the definition of rule, we know that there are $\beta_i \in [0,1]$, $i \in N$, such that $R_i = \beta_ic_i, i \in N$, and $\sum_{i \in N}\beta_i c_i=e$.

Consider the following chain of problems:
$$
(\mathcal{I},N,E,c,\alpha) \longrightarrow (\mathcal{I},N(q),E,c(q),\alpha)
$$
where the first problem is the original, the second is the problem in which each claimant $i$ is split into a number of identical claimants $k_i$, $k_i \in \mathbb{N}_+$, with claims exactly equal to $q \in \mathbb{R}_+$. We now distinguish two cases:

\begin{enumerate}
\item $c_1,c_2,\ldots,c_n \in \mathbb{Q}_+$. In this case, there exists $q \in \mathbb{Q}_+$ such that $c_i = k_iq, k_i \in \mathbb{N}_+, i \in N$. Now, by $PEFF$ and $ETE$, we have that 
$$
R_j(\mathcal{I},N(q),E,c(q),\alpha)=\beta q, j \in N(q).
$$ 
On the other hand, by $NMS$, it holds for every $i \in N$ that
$$
\beta_i c_i = k_i \beta q = \beta c_i \Rightarrow \beta_i = \beta
$$

\item $c_1,c_2,\ldots,c_n \in \mathbb{R}_+$. In this case, for each $\varepsilon >0$, there exists $q \in \mathbb{R}_+$ such that $c_i = k_i(q)q + \varepsilon_i(q), k_i(q) \in \mathbb{N}_+$, and $\varepsilon_i(q) < \frac{\varepsilon}{n}$, for all $i \in N$.

Now, by $ETE$,  we have the following equality for the second problem:
$$
\left(\sum_{i \in N} k_i(q)\right)\beta(q) q + \sum_{i \in N} \delta_i(q)\varepsilon_i(q)=e,
$$
where $\beta(q)\in [0,1]$, and $\delta_i(q) \in [0,1]$ for all $i \in N$. This equality can be written as follows:
$$
\beta(q)\sum_{i \in N} \left(c_i - \varepsilon_i(q)\right) + \sum_{i \in N} \delta_i(q)\varepsilon_i(q)=e,
$$
or equivalently,
$$
\frac{e}{\sum_{i \in N} c_i}-\beta(q) = \frac{\sum_{i \in N}\left\{ (\delta_i(q) - \beta(q))\varepsilon_i(q)\right\}}{\sum_{i \in N} c_i},
$$
taking limits on both sides when $q$ goes to zero, we obtain that $\lim_{q \rightarrow 0^+}\beta(q) = \frac{e}{\sum_{i \in N} c_i}$.

On the other hand, by $NMS$, for each $q$ and for each $i \in N$,
$$
\beta_i c_i = k_i(q) \beta(q) q + \delta_i(q)\varepsilon_i(q) = \beta(q) c_i + (\delta_i(q) - \beta(q))\varepsilon_i(q).
$$

Since $\lim_{q \rightarrow 0^+}\beta(q) = \frac{e}{\sum_{i \in N} c_i}$, $\beta_i = \frac{e}{\sum_{i \in N} c_i}$, for each $i \in N$.
\end{enumerate}
\end{proof}
\end{lemma}

\begin{lemma}\label{lemma2}
For each problem $(\mathcal{I},N,E,c,\alpha) \in \mathcal{MAC}$, and each rule $R$ that satisfies $PEFF$, $ETE$, $GMA$ and $NMS$, if for each $N' \subset N$ with $|N'|=|N|-1$, we have $R_i(\mathcal{I},N,E,c,\alpha) =CPA_i(\mathcal{I}',N',E'^R,c|_{N'},\alpha)$ for all $i \in N'$, then $R(\mathcal{I},N,E,c,\alpha) = CPA(\mathcal{I},N,E,c,\alpha)$.
\begin{proof}
We first prove that if there is $R_i=R_i(\mathcal{I},N,E,c,\alpha) =CPA_i(\mathcal{I},N,E,c,\alpha)$, then the result holds. Indeed, let us consider $R$ in the conditions of the statement, and $R_i = CPA_i(\mathcal{I},N,E,c,\alpha)$. We now consider $N' = N \backslash\{i \}$, since $R_i = CPA_i(\mathcal{I},N,E,c,\alpha)$, 
$$
(\mathcal{I}',N',E'^R,c|_{N'},\alpha) =(\mathcal{I}',N',E'^{CPA},c|_{N'},\alpha).
$$
By hypothesis, we have that for all $k \in N'$,  
$$
R_k(\mathcal{I},N,E,c,\alpha) =CPA_k(\mathcal{I}',N',E'^R,c|_{N'},\alpha)=CPA_k(\mathcal{I}',N',E'^{CPA},c|_{N'},\alpha).
$$
Moreover, since CPA satisfies consistency,
$$
CPA_k(\mathcal{I},N,E,c,\alpha) = CPA_k(\mathcal{I}',N',E'^{CPA},c|_{N'},\alpha) \text{ for all } k \in N'.
$$
Therefore, $CPA_k(\mathcal{I},N,E,c,\alpha)=R_k(\mathcal{I},N,E,c,\alpha)$ for all $k \in N'$.

Let us consider $R$ in the conditions of the statement and we assume without loss of generality that $\beta_1=\frac{R_1}{c_1} \leq \beta_2=\frac{R_2}{c_2} \leq \ldots \leq \beta_{|N|}=\frac{R_{|N|}}{c_{|N|}}$, where for the sake of simplicicty we denote $R_k(\mathcal{I},N,E,c,\alpha)$ by $R_k$ for each $k \in N$.

\medskip
First, for $\alpha(1)$, for every $i \in \alpha(1)$, we take $\gamma_i >0$ such that $\gamma_i\sum_{j \in N:i \in \alpha(j)} c_j = e_i$. We now define $\gamma_1 = \min\{\gamma_i: i \in \alpha(1)\}$, and we assume that $\gamma_1$ is without loss of generality obtained for issue 1.

\medskip
Second, $\beta_1 \leq \gamma_1$, otherwise we would have that 
$$
\sum_{j:1 \in \alpha(j)}\beta_jc_j \geq \beta_1 \sum_{j:1 \in \alpha(j)}c_j >\gamma_1 \sum_{j:1 \in \alpha(j)}c_j =e_1,
$$
which is a contradiction.

\medskip
Third, by Lemma \ref{lemma} and $GMA$, $R_1 \geq \min\{\frac{c_1}{\sum_{j:i\in \alpha(j)}c_j}e_i: i \in \alpha(1)\}=\gamma_1 c_1$. Therefore, $\beta_1 = \gamma_1$. Now, by $PEFF$, $\beta_j = \gamma_1$ for all $j \in N$ such that $1 \in \alpha(j)$.

\medskip
Fourth, for each $N' \subset N$ with $1 \in N'$ and  $|N'|=|N|-1$, by hypothesis and the definition of CPA, we have that $\beta_1$ coincides with the $\lambda^1$'s of the iterative procedures for calculating each $CPA(\mathcal{I}',N',E'^R,c|_{N'},\alpha)$. This implies that $\beta_1 =\lambda^1= \min\{\lambda^1_{i} :i \in \mathcal{I}'^1\}$, for each $N' = N \backslash\{k\}, k \in N\backslash \{1\}$, where
$$
\lambda^1_{i}\sum_{j \in \mathcal{N}'^1:i \in \alpha(j)}c_j^1 = e_i^1 - \delta(i,k)\beta_k c_k, \forall i \in \mathcal{I}',
$$
where $\delta(i,k) = 1$ if $i \in \alpha(k)$, and $0$ otherwise. Since $R_1 = \beta_1c_1=CPA_i(\mathcal{I}',N',E'^R,c|_{N'},\alpha)$, $\arg\min\{\lambda^1_{i} :i \in \mathcal{I}'^1\} \in \alpha(1)$, otherwise claimant 1 would obtain more than $\beta_1c_1$ which is a contradiction. In particular, the minimum will be attained, although possibly among others, at issue 1, because $\beta_j = \gamma_1$ for all $j \in N$ such that $1 \in \alpha(j)$.

\medskip
On the other hand, by definition of $CPA$, we have that $\lambda^1= \min\{\lambda^1_{i} :i \in \mathcal{I}^1\}$, where 
$$
\lambda^1_{i}\sum_{j \in \mathcal{N}^1:i \in \alpha(j)}c_j^1 = e_i^1, \forall i \in \mathcal{I}.
$$
This $\lambda^1$ can be also obtain by solving the following simple linear program:
$$
\begin{array}{rll}
\lambda^1= & \max & \lambda\\
                  & \text{subject to} & \lambda\sum_{j \in \mathcal{N}^1:i \in \alpha(j)}c_j^1 \leq e_i^1, \forall i \in \mathcal{I}\\
                  &  & \lambda \geq 0 
\end{array}
$$

It is obvious that $\lambda^1 \leq \gamma_1$, because of the definition of $\gamma_1$. Now, since $R$ satisfies $PEFF$, $\gamma_1$ is a feasible solution for the linear program above and $\lambda^1 \leq \gamma_1$, $\gamma_1$ is an optimal solution of the problem. Therefore, the inequality associated with issue 1 is saturated in the optimal solution and by definition of $CPA$ claimant 1 will obtain $\gamma_1c_1 = \beta_1c_1=R_1$, i.e., $R_1(\mathcal{I},N,E,c,\alpha) = CPA_1(\mathcal{I},N,E,c,\alpha)$.
\end{proof}
\end{lemma}

\begin{theorem}\label{char1}
CPA is the only rule that satisfies $PEFF$, $ETE$, $GMA$, $CONS$, and $NMS$.
\begin{proof}
We distinguish three cases, depending on the number of claimants in the problem.
\begin{enumerate}
\item $|N|=1$. In this case, all rules that satisfy $PEFF$ coincide with CPA.
\item $|N|=2$. We distinguish two cases:
\begin{enumerate}
\item $\alpha(1) \cap \alpha(2) = \varnothing$. In this situation, since the rule satisfies $PEFF$ we can consider two separate problems of only one claimant each. Now by applying the case $|N|=1$, all rules that satisfy $PEFF$ coincide with CPA.
\item $\alpha(1) \cap \alpha(2) \neq \varnothing$. We consider other two cases:
\begin{enumerate}
\item $\alpha(1) = \alpha(2)$. By $GMA$, Lemma \ref{lemma}, and the definition of rule,
$$
R_1 = \frac{c_1}{c_1+c_2}e, \text{ and } R_2 = \frac{c_2}{c_1+c_2}e,
$$
where $e$ is the minimum of the available amounts of the issues. 

\item $\alpha(1) \neq \alpha(2)$. First note that by Lemma \ref{lemma}, we know that for every single issue we obtain the proportional distribution of the available amount among the corresponding claimants. Therefore, in order to apply $GMA$, we can consider without loss of generality the two situations shown in Figure \ref{2claimants}.
\begin{figure}
\begin{center}
\begin{tikzpicture}
\node[draw] (c1) at (1.5,0) {1};
\node[draw] (c2) at (3.5,0) {2};

\node[draw] (e1) at (0.5,2.5) {$e_1$};
\node[draw] (e2) at (2.5,2.5) {$e_2$};
\node[draw] (e3) at (4.5,2.5) {$e_3$};

\draw [->] (c1) -- (e1);
\draw [->] (c1) -- (e2);
\draw [->] (c2) -- (e2);
\draw [->] (c2) -- (e3);

\node (d1) at (2.5,-1) {(A)};

\node[draw] (a1) at (8.5,0) {1};
\node[draw] (a2) at (10.5,0) {2};

\node[draw] (b1) at (8.5,2.5) {$e_1$};
\node[draw] (b2) at (10.5,2.5) {$e_2$};

\draw [->] (a1) -- (b1);
\draw [->] (a2) -- (b1);
\draw [->] (a2) -- (b2);

\node (d2) at (9.5,-1) {(B)};

\end{tikzpicture}
\end{center}
\caption{Basic 2-claimants situations when $\alpha(1) \cap \alpha(2) \neq \varnothing$ and $\alpha(1) \neq \alpha(2)$.}\label{2claimants}
\end{figure}
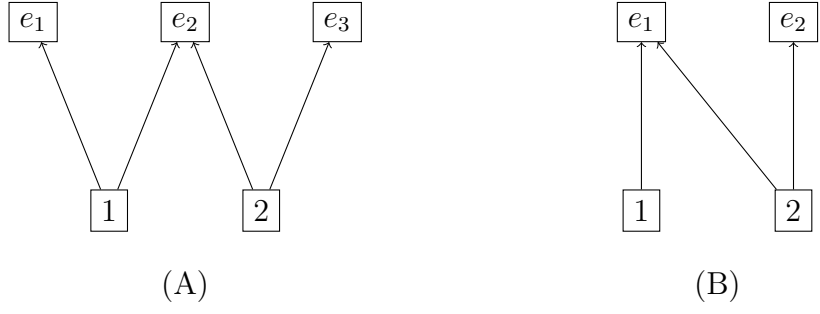

\medskip
We next analyze the two situations in Figure \ref{2claimants}:
\medskip
\begin{enumerate}
\item By $GMA$, $c_1 \geq e_1$, $c_2 \geq e_3$, and $c_1+c_2 \geq e_2$, 
$$
R_1 \geq \min\{e_1,\frac{c_1}{c_1+c_2}e_2\}, \text{ and } R_2 \geq \min\{\frac{c_2}{c_1+c_2}e_2, e_3\}
$$
If $\min\{e_1,\frac{c_1}{c_1+c_2}e_2\}=e_1$, then $R_1=e_1$, and by $PEFF$ $R_2 = \min\{c_2,e_2-e_1,e_3\}$. If $\min\{e_1,\frac{c_1}{c_1+c_2}e_2\}=\frac{c_1}{c_1+c_2}e_2$, then we have two possibilities:
\medskip
\begin{itemize}
\item $\min\{\frac{c_2}{c_1+c_2}e_2, e_3\}=e_3$, then $R_2=e_3$, and by $PEFF$ $R_1 = \min\{c_1,e_2-e_3,e_1\}$.
\medskip
\item $\min\{\frac{c_2}{c_1+c_2}e_2, e_3\}=\frac{c_2}{c_1+c_2}e_2$, then $R_1=\frac{c_1}{c_1+c_2}e_2$ and $R_2=\frac{c_2}{c_1+c_2}e_2$.
\end{itemize}
\medskip
\item By $GMA$, $c_1 \geq e_1$, $c_2 \geq e_3$, and $c_1+c_2 \geq e_2$, 
$$
R_1 \geq \frac{c_1}{c_1+c_2}e_1, \text{ and } R_2 \geq \min\{\frac{c_2}{c_1+c_2}e_1, e_2\}
$$
Now reasoning as in the previous case,
$$
R_1 = \frac{c_1}{c_1+c_2}e_1, \text{ and } R_2 = \frac{c_2}{c_1+c_2}e_1
$$
or
$$
R_1 = \min\{c_1,e_1-e_2\}, \text{ and } R_2 = e_2
$$
\end{enumerate}
\end{enumerate}
\end{enumerate}
\item $|N|= 3$. Let $R$ be a rule that satisfies $PEFF$, $ETE$, $GMA$, $CONS$, and $NMS$, and let $(\mathcal{I},N,E,c,\alpha) \in \mathcal{MAC}$, then we have that
$$
R(\mathcal{I},N,E,c,\alpha) = CPA(\mathcal{I},N,E,c,\alpha).
$$

Indeed, for each $N'=\{i_1,i_2\} \subset N$ such that $|N'|=2$, since $R$ satisfies CONS, 
$$
R_{i_k}(\mathcal{I}',N',E'^R,c|_{N'},\alpha) = R_{i_k}(\mathcal{I},N,E,c,\alpha), k =1,2,
$$
and since $|N'|=2$, we have that
$$
R_{i_k}(\mathcal{I}',N',E'^R,c|_{N'},\alpha) = CPA_{i_k}(\mathcal{I}',N',E'^R,c|_{N'},\alpha), k =1,2.
$$
Since we can take all possible $N'=\{i_1,i_2\} \subset N$, by Lemma \ref{lemma2}
$$
R(\mathcal{I},N,E,c,\alpha)=CPA(\mathcal{I},N,E,c,\alpha).
$$
\item $|N|\leq k$. Let us suppose that for each $(\mathcal{I},N,E,c,\alpha)$ with $|N| \leq k$, $R(\mathcal{I},N,E,c,\alpha)=CPA(\mathcal{I},N,E,c,\alpha)$.
\item $|N|= k+1$. For each $N' \subset N$ such that $|N'|=k$, since $R$ satisfies CONS,
$$
R_{i}(\mathcal{I}',N',E'^R,c|_{N'},\alpha) = R_{i}(\mathcal{I},N,E,c,\alpha), i \in N'.
$$
and since $|N'|\leq k$, we have that
$$
R_{i}(\mathcal{I}',N',E'^R,c|_{N'},\alpha) = CPA_{i}(\mathcal{I}',N',E'^R,c|_{N'},\alpha), i \in N'.
$$
Finally, since we can take all possible $N'\subset N$ with $|N'|=k$, by Lemma \ref{lemma2},
$$
R(\mathcal{I},N,E,c,\alpha)=CPA(\mathcal{I},N,E,c,\alpha).
$$
\end{enumerate}
\end{proof}
\end{theorem}

\begin{proposition}
Properties in Theorem \ref{char1} are logically independent.
\begin{proof}We consider the four posibilities:
\begin{itemize}
\item (No $PEFF$) The null rule satisfies all properties but $PEFF$.
\item (No $ETE$) Consider an order on the set of claimants and a rule which reimburses each claimant all that can be, in that order, until it is not possible to do it. If we assume that when a claimant splits into several new claimants or some claimants leave, the order in which the claims are attended is preserved, then this rule satisfies $PEFF$, $GMA$, $CONS$, and $NMS$, but not $ETE$.
\item (No $GMA$) Consider a rule that has two phases. In the first phase, each issue is distributed proportionally, but only among those claimants that only demand the corresponding issue. In the second phase, the amounts of each issue are updated down accordingly, and distributed among the rest of the claimants applying CPA. This rule satisfies $PEFF$, $ETE$, $CONS$, and $NMS$, but not $GMA$.
\item (No $CONS$) For every problem $(\mathcal{I},N,E,c,\alpha)\in \mathcal{MAC}$, consider the following rule defined in two steps:
\begin{enumerate}
\item First, we allocate to each claimant $j$, $\min\left\{R_j\left(\{i\},N_i, e_i, c|_{N_i}\right) : i \in \alpha(j)\right\}$.
\item Next, we revise down the available amounts of issues and the claims, and we assume without loss of generality that $e'_1 \leq e'_2\leq \ldots \leq e'_m$. Then we begin to distribute each state proportionally among the claimants, starting from the smallest to the largest quantity available. It is not until one state has been fully distributed or the claimants fully satisfied that we move on to the next updating the claims. We continue until all the states have been distributed as much as possible.
\end{enumerate}
The allocation to each claimant is the sum of everything that she has obtained in each of the steps of the procedure described.

By definition this rule satisfies $GMA$, $PEFF$, and $ETE$. Moreover, using arguments similar to those used in Theorem \ref{CPAprop}, it can be shown that this rule satisfies $NMS$. However, it does not satisfies $CONS$ since this rule does not coincide with CPA, and if we consider reduced problems with $|N'|=2$, by $PEFF$, $ETE$, $GMA$, and $NMS$, we obtain the allocations prescribed by CPA. 
\item (No $NMS$) The CEA rule for MAC satisfies all properties but $NMS$ (Acosta-Vega et al., 2021).
\end{itemize}
\end{proof}
\end{proposition}


\section{Conclusions}
\label{conc}

In allocation problems the concept of proportionality is put into practice with the well-known proportional rule. This rule has been extensively studied in the literature from many different point of views and for many allocation models. Focusing on bankruptcy models and their extensions to the multi-issue case, the proportional rule has been characterized in the context of bankruptcy problems in Chun (1988) and 
de Frutos (1999). In both papers, non-manipulability plays a central role in the axiomatic characterization of the proportional rule. For multi-issue allocation problems, Ju et al. (2007) and Moreno-Ternero (2009) introduce two different definitions of proportional rule following two different approaches. Moreover, Ju et al. (2007) and Berganti\~nos et al. (2010) provides charaterizations of both proportional rules. Again, in both approaches, non-manipulability is an essential property. In this paper, we introduce a definition of proportional rule, that we call constrained proportional awards rule, for multi-issue allocation problems with crossed claims and provide a characterization of it. Once again, non-manipulability is used. Therefore, we fill a gap in the literature of proportional distributions in allocation problems in line with the previous studies. 

Futher research will include the introduction and analysis, in this context, of the ramdon arrival rule (O'Neill, 1982). This rule is related to the well-known Shapley value (Shapley, 1953), see Algaba et al. (2019b) for an updating on theoretical and applied aspects about this value. The stydy of the Talmud rule introduced by Aumann and Maschler (1985) would be also of interest in this setting.

\section*{Acknowledgment}
This work is part of the R\&D\&I project grant PGC2018-097965-B-I00, funded by MCIN/ AEI/10.13039/501100011033/ and by "ERDF A way of making Europe"/EU. The authors are grateful for this financial support. Encarnaci\'on Algaba also acknowledges financial suport from IDEXLYON from Universit\'e de Lyon (project INDEPTH) within the Program Investissements d’Avenir (ANR-16-IDEX-0005). Joaqu\'in S\'anchez-Soriano also acknowledges financial suport from the Generalitat Valenciana under the project PROMETEO/2021/063.

\end{document}